\documentclass{article}

\usepackage{graphicx}

\newcommand{\beq}{\begin{equation}}
\newcommand{\eeq}{\end{equation}}

\begin{document}

\title{On the real zeroes of the Hurwitz
zeta-function and Bernoulli polynomials}

\maketitle

\begin{center}

{\bf A.P.Veselov $^{*,**}$ and J.P.Ward $^{*}$}

\bigskip

{\it $^{*}$Department of Mathematical Sciences, Loughborough University,
Loughborough, Leicestershire, LE 11 3TU, UK
}

\bigskip

{\it $^{**}$Landau Institute for Theoretical Physics, Kosygina 2,

    Moscow, 117940,  Russia

\bigskip

e-mails:  A.P.Veselov@lboro.ac.uk, J.P.Ward@lboro.ac.uk
}

\end{center}

\bigskip
\noindent
{\small  {\bf Abstract.} The behaviour of real zeroes of the Hurwitz
zeta function
$$\zeta (s,a)=\sum_{r=0}^{\infty}(a+r)^{-s}\qquad\qquad a > 0$$ is
investigated. It is shown that
$\zeta (s,a)$ has no real zeroes $(s=\sigma,a)$ in the region
$a >\frac{-\sigma}{2\pi e}+\frac{1}{4\pi e}\log (-\sigma)
+1$ for large negative $\sigma$. In the region $0 < a <
\frac{-\sigma}{2\pi e}$ the zeroes are asymptotically
located at the lines $\sigma + 4a + 2m =0$ with integer $m$.
If $N(p)$ is the number of real zeroes of $\zeta(-p,a)$ with given $p$ then
$$\lim_{p\to\infty}\frac{N(p)}{p}=\frac{1}{\pi e}.$$
As a corollary we have a simple proof of Inkeri's result that the
number of real roots
of the classical Bernoulli
polynomials $B_n(x)$ for large $n$ is asymptotically equal to
$\frac{2n}{\pi e}$.

\small

\section{Introduction.}

The classical {\it Hurwitz zeta-function} is defined for any positive real $a$
as an analytic continuation of the series
$$\zeta (s,a)=\sum_{r=0}^{\infty}(a+r)^{-s}.$$
When $a=1$ it reduces to the  Riemann zeta-function.

We should note that sometimes the definition of the Hurwitz
zeta-function is restricted to
$0<a\leq 1$ (see e.g. \cite{WW,Er}). From our point of view this is
not natural and we follow the definition
of the Hurwitz zeta-function from \cite{AAR} where all positive $a$ are
allowed (cf. also the original Hurwitz paper \cite {Hur}).

It is known (see e.g. \cite{Er}, volume 1, page 27) that in the
special cases when $s$ is a negative integer this function (as a
function of the parameter
$a$) reduces, up to a factor,
to a Bernoulli polynomial: explicitly when $s=-m, m=0,1,2,3,....$

$$\zeta (-m,a)=-\frac{B_{m+1}(a)}{m+1}.$$

The Bernoulli polynomials $B_k(a)$ can be defined
through the generating function:
$$\frac{ze^{za}}{e^z-1}=\sum_{k=0}^{\infty}\frac{  B_{k}(a)}{k!}z^k$$
giving, for example:$$  B_0(a)=1\qquad   B_1(a)=a-\frac{1}{2}\qquad
B_2(a)=a^2-a+\frac{1}{6}\qquad
B_3(a)=a^3-\frac{3a^2}{2}+\frac{a}{2},...$$

Bernoulli polynomials possess many interesting
properties and arise in
many areas of mathematics (see \cite{WW,Er}).

Inkeri \cite{Inkeri} proved a remarkable fact that the number $N(n)$
of real roots of
the Bernoulli polynomials $B_n$ for large $n$ asymptotically equals
to $\frac{2n}{\pi e}$.
More precise estimates for $N(n)$ have been found by Delange
\cite{Delange1,Delange2}.

In this paper we investigate the behaviour of the real zeroes of the
Hurwitz zeta function
$\zeta(\sigma,a)$ in the upper half-plane $a>0.$
As a corollary we have a simple proof of Inkeri's result.
Our approach is different from \cite{Inkeri,Delange1,Delange2} and we
believe is more elementary.
It is based on the remarkable Hurwitz representation of the
$\zeta (s,a)$ on the interval
$0<a\le 1$ and Re$(s)=\sigma <0$:

$$\zeta(s,a)=\frac{2\Gamma (1-s)}{(2\pi
)^{1-s}}\sum_{r=1}^{\infty}\frac{\sin (2\pi ra+\frac{1}{2}\pi
s)}{r^{1-s}},$$
where $\Gamma$ is the  Euler gamma-function (see e.g \cite{Er}, volume
1, page 26).

Our main observation is that for a large negative $\sigma$ this formula
gives a good approximation for
$\zeta(\sigma,a)$ on a much wider interval: $0< a < -\frac{\sigma}{2\pi
e}$. As a result we prove that

$$\frac{\zeta(\sigma,a)}{Q(\sigma)} = \sin (2\pi a + \frac{1}{2} \pi
\sigma) + o(1) \qquad {\rm as}\quad \sigma \to
-\infty\qquad {\rm provided} \quad 0<a<-\frac{\sigma}{2\pi e}$$
where $Q(s) = \frac{2\Gamma
(1-s)}{(2\pi )^{1-s}}$.
We show also that in the region $a >-\frac{\sigma}{2\pi e}+
\frac{1}{4\pi e}\log(-\sigma)+1$ the Hurwitz zeta-function has
no real zeroes for large negative $\sigma$.

$$
\includegraphics{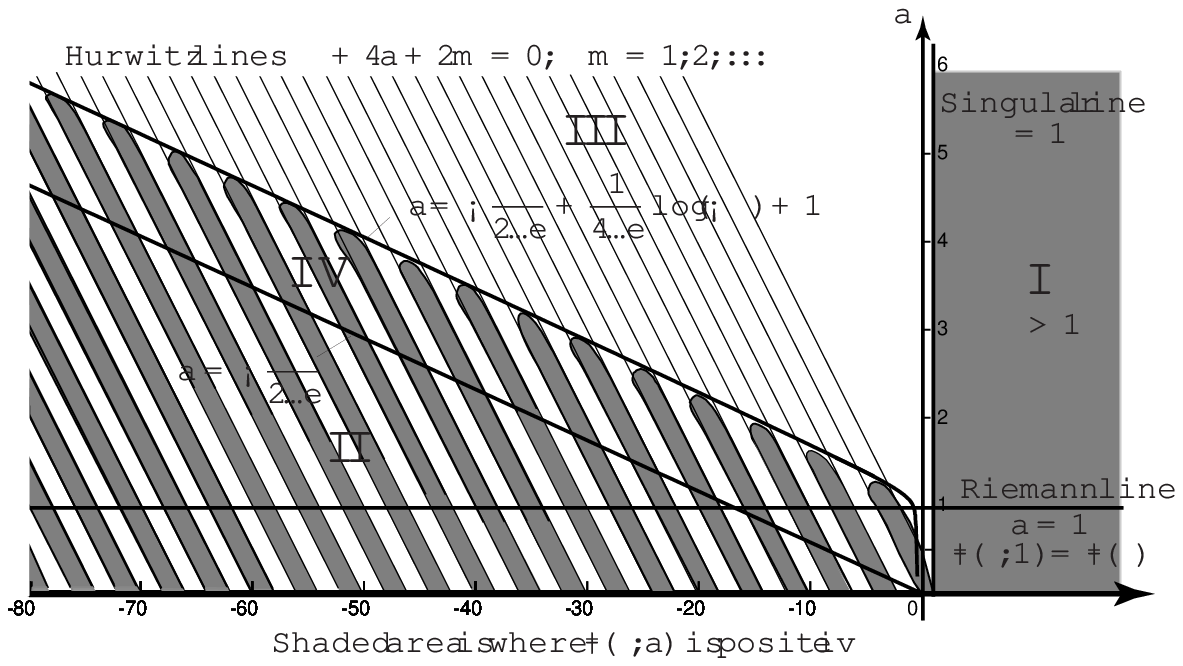}
$$
\centerline{Figure 1}

Our results are illustrated in Figure 1 which shows the
behaviour of the  real zeroes of $\zeta(\sigma,a)$. Since
$\frac{1}{2\pi e}$ is a small number we have chosen different
scales on the axes to make the picture more illuminating. Notice that
   $\zeta(s,1)$ coincides with the
  Riemann zeta-function $\zeta(s)$ which only has real zeroes
(for $s<0$) at
negative even numbers $s = -2, -4, -6,...$. Also, using the
well-known identity $\zeta(s,1/2) =
(2^s-1)\zeta(s)$, we see that the only real zeroes of $\zeta (s,a)$
on $a=1$ are $s = -2, -4, -6,...$ and on
$a=\frac{1}{2}$ they are $s = 0, -2, -4,...$ As we have mentioned
above when $s$ is a non-positive integer
$\zeta(s,a)$ reduces to a Bernoulli polynomial. We have used
this fact to compute numerically the corresponding values of
$\zeta(s,a)$ and to draw the picture in the regions II and
IV. When
$s=1$  the Hurwitz zeta-function has a simple pole with the residue 1
and for $\sigma = {\rm Re}(s) >1$ it is given by the
convergent series with positive elements and therefore has no zeroes.

\section{Asymptotic behaviour of the Hurwitz zeta-function
$\zeta(\sigma,a)$ for  large negative $\sigma$.}
\vskip0.2in

The  Hurwitz zeta-function (or generalised Riemann zeta-function) is
defined as a
series

$$\zeta (s,a)=\sum_{r=0}^{\infty}(a+r)^{-s}\qquad\qquad a>0$$
in the complex domain Re$(s) > 1$ and can be analytically continued to a
meromorphic function in the
whole complex plane with the only pole at $s = 1$ (see
\cite{WW},\cite{Er},\cite{AAR}).
When $a=1$ it reduces to the  Riemann zeta-function
$$\zeta (s)=\sum_{k=1}^{\infty} k^{-s}.$$

The Hurwitz zeta-function can be extended to the whole of the complex
$s$-plane through the formula
$$\zeta (s,a)=-\frac{\Gamma (1-s)}{2\pi {\rm
i}}\int_{\infty}^{(0+)}\frac{(-z)^{s-1}e^{-az}}{1-e^{-z}}{\rm d}z\qquad
a>0$$ in which the integral is taken over a curve starting at \lq
infinity\rq\ on the real axis,
encircles the origin in a positive direction and returns to the starting
point (see \cite{Er}). By using an
alternative integral formulation for $\zeta (s,a)$ it can be shown that
$\zeta (s,a)$ is analytic everywhere
except for the simple pole at $s=1$.

The Hurwitz zeta-function obviously satisfies the functional relation:
\beq
\zeta (s,a)=\zeta (s,n+a)+\sum_{r=0}^{n-1}(r+a)^{-s}\qquad\qquad n=
1,2,\dots
\label{Feqn}
\eeq
    Since each term in this
relation is analytic we can  assume this relation is true for the whole of
the complex
$s$ plane, except for $s=1$.

In this paper we restrict ourselves to the case when $s$ is real: $s=
\sigma \in \bf{R}$.
When $\sigma$ is negative Hurwitz has found the following Fourier
representation for $\zeta(\sigma,a)$
on the interval $0 < a\le 1$:

\beq
\label{H}
\zeta(\sigma,a)=\frac{2\Gamma (1-\sigma)}{(2\pi
)^{1-\sigma}}\sum_{r=1}^{\infty}\frac{\sin (2\pi
ra+\frac{1}{2}\pi \sigma)}{r^{1-\sigma}}
\eeq
   From this formula we see that
$$\frac{\zeta(\sigma,a)}{Q(\sigma)} = \sin (2\pi a+\frac{1}{2}\pi
\sigma) + o(1)
\qquad\qquad {\rm when}\quad \sigma \to
-\infty,$$
where $Q(\sigma) = \frac{2\Gamma (1-\sigma)}{(2\pi )^{1-\sigma}}$.
Our first theorem proves that this is actually true on a much larger interval.
%\vskip0.1in\noindent

As part of the theorem proofs we will use the following inequality
for the function
$S(p,n)=1^p+2^p+\dots +n^p$:\quad
$$S(p,n)<n^p\left(\frac{1-e^{-p}}{1-e^{-p/n}}\right)$$
Indeed,

$$S(p ,n)=n^{p }\left(1+(1-\frac{1}{n})^{p }+(1-\frac{2}{n})^{p }+\dots
+(1-\frac{n-1}{n})^{p }\right)$$ But since
$1-x<e^{-x}$ for $x<1$
we have
$$1-\frac{1}{n}<e^{-1/n}\qquad 1-\frac{2}{n}<e^{-2/n}\qquad\dots\qquad
1-\frac{n-1}{n}<e^{-(n-1)/n}$$
and therefore
$$S(p ,n)<n^{p }\left(1+e^{-p /n}+e^{-2p /n}+\dots +e^{-(n-1)p /n}\right)
= n^{p }\left(\frac{1-e^{-p }}{1-e^{-p/n}}\right).$$
We shall also need an estimate for the sum of the series:
$\sum_{r=2}^{\infty}\frac{1}{r^{1+p}}$. Now
$$\frac{1}{2^{1+p}}+\frac{1}{3^{1+p}}+\frac{1}{4^{1+p}}+\dots
=\frac{1}{2^{p}}\left(\frac{1}{2}+\frac{1}{3}(\frac{2}{3})^p+\frac{1}{4}(\frac{2}{4})^p+\dots\right)$$
But $(\frac{2}{n})^p<\frac{1}{n}$ if $p>4$ and $n\ge 3$ so we easily deduce
$$\sum_{r=2}^{\infty}\frac{1}{r^{1+p}}<(\zeta
(2)-\frac{3}{4})2^{-p}\qquad p>4$$

\noindent Finally, we will also make particular use of Stirling's
inequality for the gamma function:
$$(2\pi p)^{\frac{1}{2}}p^pe^{-p}<\Gamma (1+p)<(2\pi
p)^{\frac{1}{2}}p^pe^{-p}e^{\frac{1}{12p}}\qquad p>\!\!>1$$
{\bf Theorem 1.} {\it Let $\sigma = -p, p \ge 0$ and $0 < a < \alpha
p$ for some positive $\alpha$ then
the Hurwitz zeta-function satisfies the inequality}
\beq
\label{Ineq}
\left|\frac{\zeta (-p,a)}{Q(-p)}-\sin (2\pi a-\frac{1}{2}\pi p)\right| <
C_1 p^{-1/2}(2\pi
e\alpha)^p + C_2 2^{-p},
\eeq
{\it where $C_1, C_2$ are  constants, which do not depend on $p$. In
particular,  on the interval
$0 < a < \frac{1}{2\pi e} p$ we have the asymptotic behaviour}
$$\frac{\zeta(-p,a)}{Q(-p)} = \sin (2\pi a-\frac{1}{2}\pi p) + o(1)\qquad\qquad
{\rm when}\quad p \to\infty.$$

{\bf Proof}\\

Let us represent $a$ as $n+b$,\  $\ 0<b \leq 1$ and with $n$ integer.
It follows from the functional relation (\ref{Feqn}) that
$$\left|\frac{\zeta
(\sigma,a)}{Q(\sigma)}-\frac{\zeta
(\sigma,b)}{Q(\sigma)}\right|=\left|\frac{\zeta
(\sigma,n+b)}{Q(\sigma)}-\frac{\zeta
(\sigma,b)}{Q(\sigma)}\right|\le \frac{1}{|
Q(\sigma)|}\sum_{r=0}^{n-1}(r+b)^{-\sigma}$$
Since $0<b \le 1$  we
obviously have
$$\sum_{r=0}^{n-1}(r+b)^{-\sigma }\le S(p,n)\qquad\qquad 0<p=-\sigma $$
and, as we obtained above,
$$S(p,n)<n^p\left(\frac{1-e^{-p}}{1-e^{-p/n}}\right).$$
Also, from the Stirling formula for the $\Gamma$-function we have the
following asymptotically
exact inequality
$\Gamma(1+p) > (2\pi p)^{1/2} p^p e^{-p}$
and therefore
$$\frac{1}{Q(-p)} = \frac{{(2\pi )^{1+p}}}{2\Gamma (1+p)} < \{\frac{2\pi
e}{p}\}^p\frac{\pi}{\sqrt{2\pi p}}$$
Thus for a large $p$
$$
\left|\frac{\zeta (-p,a)}{Q(-p)}-\frac{\zeta (-p,b)}{Q(-p)}\right| <
\left\{\frac{2\pi en}{p}\right\}^p\frac{\pi}{\sqrt{2\pi
p}}\left(\frac{1-e^{-p}}{1-e^{-p/n}}\right)\qquad\qquad\qquad
$$
\beq
\qquad\qquad\qquad\qquad\qquad <\left\{\frac{2\pi
en}{p}\right\}^p\frac{\pi}{\sqrt{2\pi
p}}\frac{1}{(1-e^{-\frac{1}{\alpha}})}\qquad {\rm if}\quad \frac{n}{p}<\alpha
\label{est1}
\eeq
which is true since $\frac{n}{p}<\frac{a}{p}<\alpha$ by assumption.

\noindent However, from Hurwitz' formula (\ref{H}) it follows that
$$\frac{\zeta (-p,b)}{Q(-p)} = \sum_{r=1}^{\infty}\frac{\sin (2\pi
rb-\frac{1}{2}\pi p)}{r^{1+p}} = \sin (2\pi b - \frac{1}{2}\pi p) +
\frac{\sin (4\pi
b-\frac{1}{2}\pi p)}{2^{1+p}} + \frac{\sin (6\pi b-\frac{1}{2}\pi
p)}{3^{1+p}} +....$$
Therefore
\beq
\left|\frac{\zeta (-p,b)}{Q(-p)}-\sin (2\pi b - \frac{1}{2}\pi
p)\right|=\left|\frac{\zeta (-p,b)}{Q(-p)}-\sin (2\pi a -
\frac{1}{2}\pi p)\right| < \sum_{r=2}^{\infty} r^{-p-1}
< (\zeta(2)-\frac{3}{4}) 2^{-p}
\label{est2}
\eeq
if $p>4$.
The estimates (\ref{est1}) and (\ref{est2}) imply the theorem.

\noindent By a slight modification of the previous arguments we can
prove the following result.

\noindent{\bf Theorem 2}\quad {\it In the region} $$a>\frac{p}{2\pi
e}+\frac{1}{4\pi e}\log p +1$$
$\zeta (-p,a)$ {\it is negative if} $p$ {\it is sufficiently large.}

\noindent{\bf Proof}\quad From the same functional relation (\ref{Feqn})
$$\frac{ \zeta (-p,a)}{Q(-p)} <\frac{ \zeta
(-p,b)}{Q(-p)}-\frac{(a-1)^p}{Q(-p)}. $$ Now assuming that $p$ is
sufficiently large we can use Stirling's inequality
$\Gamma (1+p)<\sqrt{2\pi p}\left(\frac{p}{e}\right)^pe^{\frac{1}{12p}}$, so
$$Q(-p)=\frac{2\Gamma (1+p)}{(2\pi )^{1+p}}<\frac{\sqrt{2\pi
p}}{\pi}\left(\frac{p}{2\pi e}\right)^pe^{\frac{1}{12p}}$$
However, we know from the Hurwitz formula, that when $p\to\infty$
$$\frac{\zeta (-p,b)}{Q(-p)} = \sin (2\pi b-\frac{1}{2}\pi p) + o(1)$$
Therefore if $\frac{(a-1)^p}{Q(-p)}>1$ and $p$ large enough then $
\zeta (-p,a)<0$. But as we have shown
$$\frac{(a-1)^p}{Q(-p)}>\left(\frac{2\pi
e(a-1)}{p}\right)^p\sqrt{\frac{\pi }{2p}}e^{-\frac{1}{12p}}$$which is
greater
than 1 if $$a-1>\frac{p}{2\pi
e}\left(\frac{2p}{\pi}\right)^{\frac{1}{2p}}e^{\frac{1}{12p^2}}=\frac{p}{2\pi
e}e^{\{\frac{(\log 2p-\log \pi )}{2p}+\frac{1}{12p^2}\}}$$
Now, using the  inequality $e^x<1+ x + x^2$  for sufficiently small
$x$ we find that
to guarantee that $\frac{(a-1)^p}{Q(-p)}>1$ it is enough to demand that
$$a>\frac{p}{2\pi e}+\frac{1}{4\pi e}\log p + 1.$$ This implies the theorem.

\section{The Real Zeroes of the Hurwitz Zeta-function and Bernoulli
Polynomials.}

Let us now fix  $\sigma = -p$ and consider $\zeta(-p,a)$ as a function of $a$.
It follows from theorem 2 that the zeroes of this function for large $p$ are
located in the interval $0<a<\frac{p}{2\pi e}+\frac{1}{4\pi e}\log p +1.$
For given $p$ let $N(p)$ be the number of real zeroes of
$\zeta(-p,a)$, and $A(p)$
be the largest of these zeroes.

\noindent{\bf Theorem 3} {\it For} $p$ {\it sufficiently
large}
\beq
\label{HNN1}
\\
\frac{p-1}{2\pi e}-\frac{1}{2}<A(p)<\frac{p}{2\pi e}+\frac{1}{4\pi e}\log p +1
\eeq
\beq
\frac{p-1}{\pi e}-1<N(p)<\frac{p-1}{\pi e}+ \frac{1}{2}\log p +2\pi e +2
\label{HNN}
\eeq
    {\it The zeroes of} $\zeta(-p,a)$ {\it on the interval}
$0<a<\frac{p-1}{2\pi e}$\quad{\it are simple and close to the half-integer
lattice:}  $a=\frac{p}{4} + \frac{l}{2}, \quad l \in {\bf Z}.$

\vskip0.1in\noindent{\bf Proof}\quad Let us introduce the function
$Z_p(a) =\frac{\zeta(-p,a)}{Q(-p)}$.
   From theorem 1 it
follows that
$$ Z_p(a)= \sin (2\pi a-\frac{1}{2}\pi p) + o(1)$$on the interval
$I_p: 0<a<\frac{p}{2\pi
e}$ when $p\to\infty$. Actually this is true also for the
$k^{{\rm th}}$ derivative of $Z_p(a)$ but on a smaller interval
$I_{p-k}$. Indeed,
from the definition of the Hurwitz zeta-function we see that
$$\frac{\partial}{\partial a} \zeta(s,a) = (-s) \zeta(s+1,a).$$ From
the property of the $\Gamma$-function
$\Gamma(p+1) = p \Gamma(p)$ it follows that $$Q(-p)= \frac{2\Gamma
(1+p)}{(2\pi )^{1+p}} = \frac{p}{2\pi} Q(-p+1).$$ Thus the derivative
of $Z_p(a)$ is equal to
$$Z_p'(a) = 2 \pi Z_{p-1}(a) = 2 \pi \sin(2\pi a-\frac{1}{2}\pi
(p-1)) + o(1) =2 \pi \cos (2\pi a-\frac{1}{2}\pi p) +
o(1)$$ on the interval $I_{p-1}.$
Similarly we have for the $k^{{\rm th}}$ derivative of $Z_p(a)$
$$Z^{(k)}_p(a) = \sin ^{(k)} (2\pi a-\frac{1}{2}\pi p) + o(1)$$ on
the interval $I_{p-k}.$

In particular, on the interval $I_{p-1}$ the function $Z_p(a)$ (and
its derivative) tend to
$\sin (2\pi a-\frac{1}{2}\pi p)$ (and its derivative)  when $p\to
\infty$ which ensures that for large $p$ all the roots
of
$\zeta(-p,a)$ on this interval are simple and located near the points
$a=\frac{p}{4}+\frac{l}{2}\quad l\in  \bf{Z} $.
This implies the last statement of the theorem and the lower
estimates of (\ref{HNN1}) and (\ref{HNN}).

The
upper estimates for
$A(p)$ follows directly from theorem 2. To prove the upper estimates
for $N(p)$ we need the following simple lemma.
\\{\bf Lemma}\quad If a function $f(x)$ (with a continuous $n^{\rm
th}$ derivative)  on some interval $(a,b)$ has the
property that the sign of the
$n^{\rm th}$ derivative is constant throughout the interval then $f$
has no more than $n$  roots on this interval.

Now we  apply this lemma to the function $Z_p(a)$ on the interval
$J_p:\ (\frac{p-1}{2\pi e},\frac{p}{2\pi
e}+\frac{1}{4\pi e}\log p +1)$ to estimate the number of roots there.
The idea of this calculation is clear from  Figure
2 (in which $\kappa\equiv\frac{1}{2\pi e}$).

$$
\includegraphics{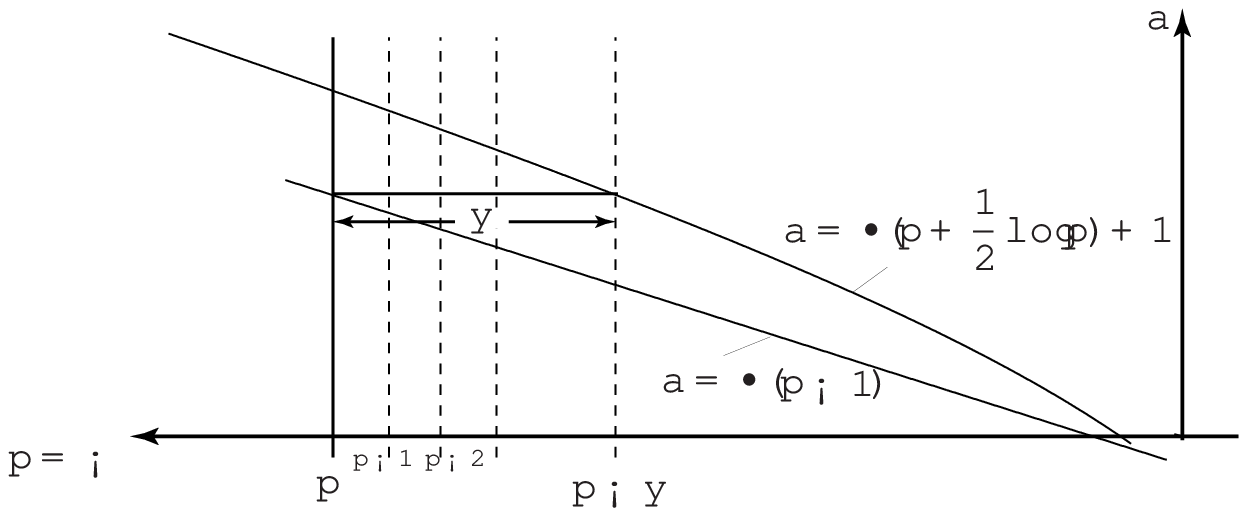}
$$
\centerline{Figure 2}

Using the fact that $$Z^{(n)}_p(a)=(2 \pi)^n Z_{p-n}(a)$$
we differentiate $Z_p(a)$ many times until we have a negative
function and then apply the lemma. As
one can see from Figure 2 if $n>[y]$ then $Z^{(n)}_p(a)$ will be
negative in the interval $J_p$ and, as such, cannot
have more than $n$ simple roots in this interval. Now
$y$ is the solution to the equation
$$ \kappa\left((p-y)+\frac{1}{2}\log (p-y)\right)+1=\kappa (p-1)$$
or
$$ -y+\frac{1}{2}\log (p-y) +(\kappa^{-1}+1) =0$$
We claim that the solution to this equation for large $p$ satisfies the
inequality
$$y<\frac{1}{2}\log p+2\pi e +1 $$
Indeed the function $F(y)= -y+\frac{1}{2}\log (p-y)+(2\pi e +1) $
is monotonically decreasing and
$$F(\frac{1}{2}\log p+2\pi e +1 )=-\frac{1}{2}\log p -2\pi e-1+\frac{1}{2}\log
(p-\frac{1}{2}\log p-2\pi e-1)+(2\pi e +1)$$
$$=\frac{1}{2}\log (1-\frac{1}{2p}\log p-\frac{2\pi e+1}{p})<0\quad\qquad$$
for  large $p.$
Thus according to the lemma $Z_p(a)$ has no more than
$\frac{1}{2}\log p+2\pi e +1$ roots on the interval
$(\frac{p-1}{2\pi e},\frac{p}{2\pi e}+\frac{1}{4\pi e}\log p +1)$.
Since on the interval
$(0, \frac{p-1}{2\pi e}]$ we have no more than $\frac{p-1}{\pi e}+1$
zeroes this
proves theorem 3.

As we have already mentioned, when $p=-\sigma =m, \ m\in \bf{Z}^+ $
the Hurwitz-zeta function reduces to certain
polynomials related in a simple way to the Bernoulli polynomials:
\beq\zeta (-m,a)=-\frac{B_{m+1}(a)}{m+1}
\label{Sym1}
\eeq
Theorem 3 applied to these special values of $p$ gives some
estimates on the real positive roots
of the Bernoulli polynomials but because of their well-known symmetry
properties:
\beq
B_m(1-a)=(-1)^mB_m(a)
\label{Sym}
\eeq
we can immediately extend this result for all real roots of $B_m(a)$.
In particular if ${\bf N}(m)$ is the number of all real roots of
$B_m(a)$ and ${\bf A}(m)$ is
the largest of these roots then from Theorem 3 it follows that for a large $m$

\beq
\label{NN1}
\\
\frac{m}{2\pi e}-\bigg(\frac{1}{\pi e}+\frac{1}{2}\bigg)<{\bf
A}(m)<\frac{m}{2\pi e}+\frac{1}{4\pi e}\log
m+\bigg(1-\frac{1}{2\pi e}\bigg)
\eeq
\beq
\frac{2m}{\pi e}-\bigg(\frac{2}{\pi e}+2\bigg)<{\bf
N}(m)<\frac{2m}{\pi e}+ \log m+\bigg(4\pi e +3-\frac{4}{\pi e}\bigg)
\label{NN}
\eeq

\noindent{\bf Corollary} ({\it K. Inkeri}).}\quad
$$\lim_{m\to\infty}\frac{{\bf N}(m)}{m}=\frac{2}{\pi e},
\quad\quad
\lim_{m\to\infty}\frac{{\bf A}(m)}{m}=\frac{1}{2\pi e}.$$

\noindent{\bf Remark.} H. Delange in \cite{Delange1, Delange2} has
found sharper
estimates for ${\bf A}(m)$ and ${\bf N}(m)$.
In particular he showed that the additional logarithmic terms exist
in both upper and lower bounds. This should be true
also for the real zeroes of the Hurwitz zeta-function but it does not
follow from our
elementary arguments.

\vskip0.2in\noindent {\bf Acknowledgements}. We are grateful to John
Gibbons and Peter
Shiu for helpful and stimulating discussions.

\end{document}